\documentclass[11pt]{article}
\usepackage{psfig,amssymb,latexsym}

\newtheorem{theorem}{Theorem}[section]

\def\slfrac#1#2{\hbox{\kern.1em %
 \raise.5ex\hbox{\the\scriptfont0 #1}\kern-.11em %
 /\kern-.15em\lower.25ex\hbox{\the\scriptfont0 #2}}}

\newcommand{\eqn}[1]{(\ref{#1})}
\newcommand{\hsp}{\hspace*{\parindent}}

\newcommand{\eeq}{\end{equation}}
\newcommand{\beql}[1]{\begin{equation}\label{#1}}
\newcommand{\bsq}{{\vrule height .9ex width .8ex depth -.1ex }}

\newcommand{\FF}{{\mathbb F}}

\newcommand{\QQ}{{\mathbb Q}}

\newcommand{\bA}{{\bf A}}

\newcommand{\bI}{{\bf I}}
\newcommand{\bR}{{\bf R}}

\makeatletter
\def\@sect#1#2#3#4#5#6[#7]#8{\ifnum #2>\c@secnumdepth
     \def\@svsec{}\else
     \refstepcounter{#1}\edef\@svsec{\csname the#1\endcsname.\hskip .75em }\fi
     \@tempskipa #5\relax
      \ifdim \@tempskipa>\z@
        \begingroup #6\relax
          \@hangfrom{\hskip #3\relax\@svsec}{\interlinepenalty \@M #8\par}%
        \endgroup
       \csname #1mark\endcsname{#7}\addcontentsline
         {toc}{#1}{\ifnum #2>\c@secnumdepth \else
                      \protect\numberline{\csname the#1\endcsname}\fi
                    #7}\else
        \def\@svsechd{#6\hskip #3\@svsec #8\csname #1mark\endcsname
                      {#7}\addcontentsline
                           {toc}{#1}{\ifnum #2>\c@secnumdepth \else
                             \protect\numberline{\csname the#1\endcsname}\fi
                       #7}}\fi
     \@xsect{#5}}
\def\@begintheorem#1#2{\it \trivlist \item[\hskip \labelsep{\bf #1\ #2.}]}

\def\plain{plain}\ifx\fmtname\plain\csname fi\endcsname
     
     \input docstrip
     \preamble

     Do not distribute the stripped version of this file.
     The checksum in the header refers to the documented version.

     \endpreamble
     \generateFile{here.sty}{t}{\from{here.doc}{}}
     \endinput
\fi
\ifcat a\noexpand @\let\next\relax\else\def\next{%
    \documentstyle[here,doc]{article}\MakePercentIgnore}\fi\next
\ifx\@Hxfloat\@Hundef\else\expandafter\endinput\fi
\let\@Hxfloat\@xfloat
\def\@xfloat#1[{\@ifnextchar{H}{\@HHfloat{#1}[}{\@Hxfloat{#1}[}}
\def\@HHfloat#1[H]{%
\expandafter\let\csname end#1\endcsname\end@Hfloat
\vskip\intextsep\vbox\bgroup\def\@captype{#1}\parindent\z@
\ignorespaces}
\def\end@Hfloat{\egroup\vskip \intextsep}

\makeatother

\catcode`\@=11
\renewcommand{\section}{
        \setcounter{equation}{0}
        \@startsection {section}{1}{\z@}{-3.5ex plus -1ex minus
        -.2ex}{2.3ex plus .2ex}{\large\bf}
        }
\catcode`@=12

\begin{document}

\begin{center}
{\Large 
{\bf A  Note on Absolute Derivations and Zeta Functions}} \\
\vspace{1.5\baselineskip}
{\em Jeffrey C. Lagarias} \\
\vspace*{.2\baselineskip}
Department of Mathematics \\
University of Michigan \\
Ann Arbor, MI 48109--1109 \\
{\tt lagarias@umich.edu}\\
\vspace*{1.5\baselineskip}

(January 17, 2005) \\
\vspace{3\baselineskip}
{\bf ABSTRACT}
\end{center}
This note answers a question raised by Kurokawa, Ochiai
and Wakayama, whether a certain operator constructed
using a notion of quantum non-commutativity of primes has
eigenvalues related to the Riemann zeta zeros.

\setlength{\baselineskip}{1.0\baselineskip}
%
%
%
%
\section{ Introduction}
\hsp
In studying the parallel between zeta functions
of number fields and
function fields over finite fields, certain properties
of number fields seem describable by 
 viewing them 
as geometric objects over the
``field with one element.'' 
Analogies in these
directions have been formalized only recently,
in Manin \cite{Ma95},  Soul\'{e} \cite{So99}, \cite{So03},
Kurokawa, Ochiai and Wakayama  \cite{KOW03}, and
Deitmar \cite{De04}. There is some earlier work,
such as Kurokawa \cite{Ku92}, which can be traced 
in the references in the papers above.

In particular, 
Kurokawa, Ochiai and Wakayama  \cite{KOW03} 
recently introduced  a notion of absolute derivation over 
the rational number  field $\QQ$. Based on this, they
proposed a measure of ``quantum non-commutativity''
of pairs of  primes over the rational field, given 
as follows. For real variables $x, y > 1$, define
\beql{100}
F(x, y) = \sum_{k=1}^{\infty} x^{k-1}\frac{y^{-x^k}}{(1 - y^{-x^k})^2}.
\eeq
Now define, for $x, y > 1$
\beql{100b}
QNC(x, y) := \frac{1}{12xy}( x(y-1)F(x,y) - y(x-1)F(y,x) ).
\eeq
The ``quantum non-commutativity'' of two primes $p$ and $q$
is defined to be $QNC(p, q)$. 
It is easy to see that $QNC(x,y) = - QNC(y, x)$, whence 
$QNC(x,x)=0$, and one has $QNC(2, 3) = 0.00220482...$,
for example.
They then raised questions (\cite[p. 580]{KOW03}) whether
there is a
connection between the quantum non-commutativity measure
and zeta functions.
Define  the  infinite skew-symmetric matrix
$\bR = \left[ \bR_{ij} \right]$ whose $(i, j)$-th entry
$$
\bR_{ij} := QNC(p_i, p_j), 
$$
where $p_i$ denotes the
$i$-th prime listed in increasing order, so that
$p_1=2, ~p_2=3,~ p_3= 5$ etc. 
In question (A) they asked 
whether it could be true (in some suitable sense)  that
\beql{101}
\det\left( \bI - \bR (s- \frac{1}{2}) \right) = c~ \xi(s),
\eeq
in which  
$\xi(s) = \frac{1}{2}s(s-1) \pi^{-\frac{s}{2}}\Gamma(\frac{s}{2})\zeta(s),$
is the Riemann $\xi$-function and $c$ is a nonzero constant.
(They proposed $c=2$.)
They also asked a more general  question (B) for 
(suitable) automorphic or Galois representations 
$\rho$, which would involve a skew-symmetric matrix $\bR(\rho)$ with
$(i,j)$-th entry 
$$
\bR_{ij}(\rho) := \frac{ \rho(p) + \rho(q)^{\ast}}{2} \bR_{ij},
$$
involving weighted version of elements 
$QNC(p_i, p_j)$, and asks whether it could be true that
\beql{101b}
\det\left( \bI - \bR(\rho) (s- \frac{1}{2}) \right) = 
c s^{m(\rho)} (s-1)^{m(\rho)}\hat{L}(s, \rho),
\eeq
where $\hat{L}(s, \rho)$ is the completed $L$-function attached to 
the representation $\rho$, and $m(\rho)$ is the multiplicity of
the trivial representation in $\rho$.

In order to make  questions (A) and (B)  well-defined one must
formulate a  suitable definition  of infinite determinant in \eqn{101}.
We take as a basic requirement of such
an infinite determinant that any  zero $s$ of
a determinant \eqn{101} must necessarily have 
$z= \frac{1}{s-\frac{1}{2}}$ belonging to 
the spectrum of $\bR$, i.e. that for this value the  
resolvent $(zI -\bR)^{-1}$ is not  a bounded operator
on the full domain of $\bR$, assumed to be a Banach space.

One consequence of this
basic requirement is that if $\bR$ acts as  a bounded operator on
some Hilbert space in \eqn{101}, then 
a positive answer to question (A) would 
necessarily imply the Riemann hypothesis.  
This follows since  $\bR$ would then be skew-adjoint, hence have
pure imaginary spectrum, whence the determinant
(assumed defined) could only vanish when $s- \frac{1}{2}$ is pure
imaginary. 
One can weaken question (A) so that it no longer
implies the Riemann hypothesis, by requiring only that 
the left side 
$\det\left( \bI - \bR (s- \frac{1}{2}) \right)$ of 
\eqn{101} 
detect all the zeta zeros that are on the
critical line $\Re(s) = \frac{1}{2}.$

This note gives a negative
answer to question (A) in both formulations.
We treat the operator $\bR$
as acting on the Hilbert space $l_2$ of column
vectors, and will show it is bounded. It 
follows that it is skew-adjoint and so has spectrum
confined to the imaginary axis. However we  show that its
spectrum cannot detect~\footnote{ If $\rho= \frac{1}{2} + i \gamma$ is a zeta
zero, the corresponding point of the spectrum  of $\bR$ is
$\lambda=-\frac{i}{\gamma}$.}
all the zeta zeros that lie on the critical line, whether
or not the Riemann hypothesis holds. 

The main point is that the quantum non-commutativity
function is so rapidly decreasing as $p, q$ increase that  
\beql{102}
\sum_{j=1}^\infty \sum_{k=1}^\infty |\bR_{jk}| < \infty,
\eeq
We show this in \S2, and deduce that 
the matrix $\bR$ defines a trace class operator on $l_2$.
The weaker condition 
\beql{102b}
\sum_{j=1}^\infty \sum_{k=1}^\infty |\bR_{jk}|^2 < \infty,
\eeq
already implies that $\bR$ is a compact operator
(in fact a Hilbert-Schmidt operator), see
\footnote{ In  Akhiezer and Glazman, the term  ``completely
continous operator'' $=$ ``compact operator''.}
Akhiezer and Glazman \cite[Sect. 28]{AG93}.
A compact operator
necessarily  has a pure discrete spectrum
with all nonzero eigenvalues of finite multiplicity,
with only limit point zero (\cite[Theorem VI.15]{RS80}. 
Since we now know $\bR$ is skew-adjoint,   
its eigenvalues, which necessarily  occur
in complex
conjugate pure imaginary pairs, and can be  
ordered by decreasing absolute value, 
$\{\pm i\lambda_j: j=1, 2, \cdots \}$.
with $\lambda_1 \ge \lambda_2 \ge \cdots >0$.
A trace class operator $\bA$ is a
compact operator with the property
that its singular values $\mu_j$
(eigenvalues of the positive self-adjoint
operator $(\bA^{\ast}\bA)^{\frac{1}{2}}$ )
satisfy 
\beql{103a} 
\sum_{j=1}^{\infty} \mu_j < \infty.
\eeq
For skew-adjoint operators $\mu_j= |\lambda_j|$,
giving the condition
\beql{103} 
\sum_{j=1}^{\infty} |\lambda_j| < \infty.
\eeq
For trace class operators $\bA$ there
 is an essentially  unique definition of
$\det(I + \bA)$
that satisfies the basic requirement, 
see B. Simon \cite{Si77},
who reviews three equivalent definitions of this determinant
(see also \cite{Si79}[Chap. 3]).
He bases his treatment on the formula  
$$
\det(I- w \bA) := \sum_{k=0}^{\infty} Tr(\bigwedge{}^k (w\bA)) = 
\sum_{k=0}^{\infty} Tr(\bigwedge{}^k \bA) w^k. 
$$
which is also presented in  Reed and Simon \cite[Sec. XIII.17, p.323]{RS78}.
This determinant is an entire function  in the variable $w$,
given by the convergent infinite product
$$
\det(I- w \bA) = \prod_j(1 - w\lambda_j(\bA)),
$$
which counts the eigenvalues $\lambda_j(\bA)$ of $\bA$ with their
algebraic multiplicity, and which 
satisfies  the basic requirement, 
Reed and Simon \cite[Theorems XIII.105(c), XIII.106]{RS78}. 
The  truth of \eqn{101} for the
trace class operator $\bR$, taking $w= s- \frac{1}{2}$,
 would imply  that  if
$ s=\frac{1}{2} + i \gamma_j$ is a zeta zero on the
critical line, 
then the two values $\lambda_j = \pm \frac{i}{\gamma_j}$ belong
to the spectrum of $\bR$.
It is well known (\cite[Chap. X]{TH86})
that a positive proportion of zeta zeros lie
on the critical line $\Re(s) = \frac{1}{2}$,
 and the  asymptotics of these zeros 
easily give
\beql{104}
\sum_{\{\gamma: \zeta(\frac{1}{2} + i \gamma)=0\}} 
\frac{1}{|\gamma|} = +\infty.
\eeq
This contradicts \eqn{103}.

In \S3 we discuss the problem of whether the notion of
``QNC'' can be modified to give a positive answer
to question (A).

%
%
%
%
\section{Proof}

Our object is to show:

\begin{theorem}~\label{th21}
The operator $\bR$ acting on the column vector space
$l_2$ defines  a trace class operator.
\end{theorem}

\paragraph{Proof.}
A bounded operator $\bA$ is trace class if $|\bA|=(\bA^{*}\bA)^{\frac{1}{2}}$
is trace class, ie. 
the positive operator $|\bA|$ has  pure discrete spectrum and the
sum of its eigenvalues converges
cf. Reed and Simon ~\cite[Sect VI.6]{RS80}. 
A necessary and sufficient condition for an operator $\bA$ to
be trace class is that for every orthonormal basis 
$\{\phi_n : 1 \le n < \infty\}$ of $l_2$  one has
\beql{202a}
\sum_{n=1}^{\infty} |\langle \bA \phi_n, \phi_n\rangle | < \infty
\eeq
see Reed and Simon \cite[Chapter VI, Ex. 26, p. 218]{RS80}.

Taking $\bA= \bR$, since it is skew-symmetric 
we have
$\bR^{\ast} \bR= - \bR^2$. 
It follows that if  $|\bR|$ is trace class,
then it has pure discrete spectrum and 
the singular values of $\bR$ are just the absolute values
of the eigenvalues of $\bR$.

We now prove \eqn{102}. We have 
$$
|QNC(p, q)| \le \frac{1}{12}(F(p,q) + F(q, p)) 
$$
Now we have  $p, q \ge 2$ so $(1 - p^{-q^k})^2 \ge \frac{9}{16}$,
whence
$$
F(p, q) \le \frac{16}{9} \sum_{k=1}^\infty p^{k-1} q^{-p^k}
\le 2 q^{-p} +   2q^{-p}(\sum_{k=2}^\infty p^{k-1} q^{p-p^k}) \le 6 q^{-p}.
$$
In the last step above we used
\footnote{Note that  $x2^{-x}$ is decreasing for 
$x \ge 2 > \frac{1}{\log 2}$.}
 (for  $k, p, q \ge 2)$
$$
p^{k-1} q^{p-p^k} \le p^{k-1}2^{-p^{k-1}} \le 2^{k-1} 2^{-2^{k-1}}
 \le 2^{-k+2}.
$$
This yields
$$
|QNC(p, q)| \le \frac{1}{2}(p^{-q}+ q^{-p}),
$$
from which  we obtain
$$
\sum_{j=1}^{\infty} \sum_{k=1}^\infty |\bR_{jk}| 
\le \sum_{m=2}^\infty \left(\sum_{n=m}^\infty m^{-n}\right) 
< \infty,
$$
as asserted.

We use \eqn{102} to verify criterion \eqn{202a}.
Let $\{e_k: 1 \le k < \infty\}$ be the standard 
orthonormal basis of column vectore of $l_2$, so that
$\bR (e_k) = \sum_{j=1}^\infty \bR_{jk}e_j$. Now let
$\phi_n = \sum_{k=1}^\infty c_{nk} e_k$
be an orthonormal basis of $l_2$, so that
$[c_{nk}]$ is a unitary matrix. Then we have 
$||\phi_n||^2 = \sum_{k=1}^\infty |c_{nk}|^2 = 1,$
and unitarity   also implies
\beql{203}
\sum_{n=1}^\infty |c_{nk}|^2 = 1.
\eeq

Now we compute 
\begin{eqnarray*}
\sum_{n=1}^{\infty} |~\langle \bR\phi_n, \phi_n\rangle ~ | &=&
\sum_{n=1}^{\infty}~|\langle \sum_{j=1}^\infty\sum_{k=1}^\infty 
c_{nk}\bR_{jk} e_j,~ \sum_{j=1}^\infty c_{nj} e_j \rangle|\\
 &\le&
\sum_{n=1}^{\infty} ~ \sum_{j=1}^\infty \sum_{k=1}^\infty
|c_{nk} \bR_{jk} \overline{c_{nj}} | \\
&\le & 
\sum_{j=1}^\infty \sum_{k=1}^\infty|\bR_{jk}|
(\left(  \sum_{n=1}^{\infty} |c_{nj}||c_{nk}|\right) \\
&\le & 
\sum_{j=1}^\infty \sum_{k=1}^\infty |\bR_{jk}|
\left( \sum_{n=1}^{\infty} \frac{1}{2}( |c_{nj}|^2 +|c_{nk}|^2)\right) \\
& \le & \sum_{j=1}^\infty \sum_{k=1}^\infty |\bR_{jk}| < \infty 
\end{eqnarray*}
as required. $~~~\bsq$

%
%
%
%
\section{Concluding Remarks}

It is a interesting question whether the concept
of  ``QNC''has a natural modification  to correct the difficulty
observed here, and possibly to  give a positive answer to question (A).
We have no proposal how to do this, but make the
following remarks.

The argument made above rests on the following
fact:  A  necessary condition 
on a  skew-symmetric compact operator $\bR$ acting on
a Hilbert space 
to have a determinant \eqn{101} satisfying 
the basic requirement that detects the zeta function zeros
is that it be a Hilbert-Schmidt operator which is not
of trace class.
In order to define \eqn{101} when $\bR$ is
a Hilbert-Schmidt operator that is not of trace class,
an extended definition of  infinite determinant is
required. There are notions of regularized 
determinant $\det(I + w\bA)$  that
apply to Hilbert-Schmidt operators $\bA$
and satisfy the basic
requirement. One such,  denoted $\det_2(I + w\bA)$,
in discussed in Simon \cite{Si77} and 
Simon \cite[Chap. 3]{Si79}.
See  Pietsch \cite[Chapters 4, 7]{Pi87} for  
further work on such  questions.

The results of Kurokawa, Ochiai and Wakayama \cite{KOW03}
were  motivated in part by the
function field case for the absolute function
field $K= \FF_q(T)$, as noted at the beginning of
their paper. We note that one might reconsider the
function field analogy, varying the base function field. 
For the (absolute) function field case $\FF_q(T)$
the corresponding matrix (and operator) $\bR \equiv 0$, but if one
allowed other function fields $K$ of genus one or
higher, then the function field analogue
of the quantity \eqn{104} also diverges.
This holds because the  function field zeta
zeros $\frac{1}{2} + i \gamma$ have $\gamma$
falling  in a finite number of arithmetic progressions
$ (\bmod~ \frac{2\pi}{\log~p})$,
so that 
$$\sum_{\gamma} \frac{1}{|\gamma|} = + \infty.$$
Thus the difficulty
above manifests itself already in the function field
case. It therefore might be useful to look 
for  formulas for quantum non-commutativity for prime
ideals in a function field $K$ of genus at least one, 
intending to construct an 
analogous matrix $\bR_{K}$. 
The operator corresponding to  $\bR_{K}$
on $l_2$ would necessarily be Hilbert-Schmidt, but  
not of trace class, if it were to have eigenvalues
$\pm \frac{i}{\gamma}$, where $\frac{1}{2} + i \gamma$
runs over 
the function field zeta zeros of $K$, counted
with multiplicity. Perhaps such study could clarify
the notion of ``QNC''.

Finally we note that if to  the sum defining the
function $F(x,y)$  in \eqn{100} the term $k=0$ were added, 
the defintion of $QNC(p, q)$ would be modified to
add the extra terms
$$
\frac{1}{12pq}(\frac{1}{q-1} - \frac{1}{p-1}).
$$ 
The resulting
modified  operator $\tilde{\bR}$ then has
$$
\sum_{i, j} |\tilde{\bR}_{ij}| = + \infty,
$$
and is a Hilbert-Schmidt
operator on $l_2$ not of trace class.

\paragraph{Acknowledgment.} The author thanks the
reviewer for helpful comments.

%
%
%
%

\end{document}